\documentclass[hypertex,12pt, reqno]{amsart}
\usepackage{fullpage}
\usepackage{hyperref}
\usepackage{graphicx}
\usepackage{amsrefs}

\theoremstyle{plain} \newtheorem*{thm*}{Theorem}
\newtheorem{thm}{Theorem}
\newtheorem{prop}[thm]{Proposition} \newtheorem{lemma}[thm]{Lemma}
\newtheorem*{prop*}{Proposition}
 \theoremstyle{definition}
\theoremstyle{remark}

\newtheorem*{dfn*}{Definition}
 \theoremstyle{remark}
 \newtheorem*{rem*}{Remark}
\newtheorem*{example*}{Example} \numberwithin{equation}{section}

  \newcommand{\bbC}{\mathbb{C}}  \newcommand{\bbN}{\mathbb{N}} \newcommand{\bbZ}{\mathbb{Z}} 
\newcommand{\f}{\mathfrak} \newcommand{\fg}{\mathfrak{g}} \newcommand{\fk}{\mathfrak{k}} \newcommand{\fp}{\mathfrak{p}} \newcommand{\fh}{\mathfrak{h}}       \newcommand{\fl}{\mathfrak{l}}  \newcommand{\fsl}{\mathfrak{sl}}  \newcommand{\fso}{\mathfrak{so}}
\newcommand{\fsp}{\mathfrak{sp}}

\newcommand{\la}{\lambda} \newcommand{\al}{\alpha} \newcommand{\be}{\beta}
\newcommand{\ga}{\gamma} \newcommand{\de}{\delta}

\newcommand{\fhd}{\mathfrak{h}^*}
\newcommand{\q}{{\bf q}}
\newcommand{\bl}{{\bf b}}
\newcommand{\ch}{\protect \mbox{ch}}

\begin{document}
\title{A generating function for Blattner's formula}

\author{Jeb F. Willenbring}
\address{
Jeb F. Willenbring \\
University of Wisconsin - Milwaukee \\
Department of Mathematical Sciences \\
P.O. Box 0413 \\
Milwaukee WI 53201-0413} \email{jw@uwm.edu}
\thanks{The first author was supported in part by NSA Grant \# H98230-05-1-0078.}

\author{Gregg J. Zuckerman}
\address{
Gregg J. Zuckerman \\
Yale University Mathematics Dept. \\
PO Box 208283 \\
New Haven, CT 06520-8283} \email{gregg@math.yale.edu}
\date{August 29, 2006}

\date{\today}
\begin{abstract}
Let $G$ be a connected, semisimple Lie group with finite center and
let $K$ be a maximal compact subgroup.  We investigate a method to
compute multiplicities of $K$-types in the discrete series using a
rational expression for a generating function obtained from
Blattner's formula.  This expression involves a product with a
character of an irreducible finite dimensional representation of $K$
and is valid for any discrete series system.  Other results include
a new proof of a symmetry of Blattner's formula, and a positivity
result for certain low rank examples.  We consider in detail the situation for $G$ of type
split $\rm G_2$. The motivation for this work came from an
attempt to understand pictures coming from Blattner's formula, some
of which we include in the paper.
\end{abstract}
\maketitle

\section{Introduction}
In \cite{HS}, a proof of a formula for the restriction of a discrete
series representation (see \cite{HC}) of a connected, linear,
semisimple Lie group to a maximal compact subgroup is given. This
formula was first conjectured by Blattner.  We recall the formula
and its context briefly, from the point of view of root system
combinatorics.

Throughout the paper, $\fg$ denotes a semisimple Lie algebra over $\bbC$ with a fixed
Cartan subalgebra $\fh$.  Let $\Phi := \Phi(\fg, \fh)$ denote the corresponding root system with Weyl group
$W_\fg$.  Choose a set,  $\Phi^+$, of positive
roots and let $\Pi := \{\al_1,\cdots,\al_r \}
\subseteq \Phi$ be the simple roots.  Let $\Phi^- = - \Phi^+$.

We assume that there exists a function $\theta:\Phi \rightarrow \bbZ_2 $
such that if $\ga_1, \ga_2 \in \Phi$ and $\ga_1 + \ga_2 \in \Phi$
then $\theta(\ga_1 + \ga_2) = \theta(\ga_1) + \theta(\ga_2)$.  This map
provides a $\bbZ_2$-gradation on $\Phi$.  We set:
\[
\begin{array}{lll}
     \Phi_c    &:=& \{ \ga \in \Phi | \theta(\ga) = 0 \}, \\
     \Phi_{nc} &:=& \{ \ga \in \Phi | \theta(\ga) = 1 \}.
\end{array}
\]
Given $\al \in \Phi$, set $ \fg_\al=\{X\in\fg \; | \; [H,X]=\al(H) X
\; \forall H \in \f h \}$.  Let $\fk := \fh \oplus \sum_{\al \in
\Phi_c} \fg_\al$ and \\ $\fp := \sum_{\al \in \Phi_{nc}} \fg_\al$.
Then, $\fk$ will be a reductive symmetric subalgebra of $\fg$ with
$\fg = \fk \oplus \fp$ the corresponding Cartan decomposition of
$\fg$.  As defined, $\fh$ is a Cartan subalgebra for $\fk$ so rank
$\fk$ = rank $\fg$. Each equal rank symmetric pair corresponds to at
least one $\bbZ_2$-gradation in this manner, and conversely.

We shall refer to the elements of $\Phi_c$ (resp. $\Phi_{nc}$) as
compact (resp. noncompact).  The compact roots are a sub-root system
of $\Phi$.  Let $\Phi^+_c := \Phi^+ \cap
\Phi_{c}$, $\Phi^+_{nc} := \Phi^+ \cap \Phi_{nc}$, $\Pi_c := \Pi
\cap \Phi_{c}$, and $\Pi_{nc} := \Pi \cap \Phi_{nc}$.
Set $\rho_\fg := \rho_c + \rho_{nc}$ where
$\rho_c := \frac{1}{2} \sum_{\al \in \Phi^+_c} \al$ and
$\rho_{nc} = \frac{1}{2} \sum_{\al \in \Phi^+_{nc}} \al$.  If there is no
subscript, we mean $\rho = \rho_c$.

We remark that the $\bbZ_2$-gradation $\theta$ is determined
by its restriction to $\Pi$.  Furthermore, to any set partition $\Pi
= \Pi_1 \biguplus \Pi_2$ there exists a unique $\bbZ_2$-gradation on
$\Phi$ such that $\Pi_c = \Pi_1$ and $\Pi_{nc} = \Pi_2$.

We denote the Killing form on $\fg$ by $(,)$, which restricts to a
nondegenerate form on $\fh$.  Using this form we may
define $\iota: \f h \rightarrow \fhd$ by $\iota(X)(-) = (X, -)$ ($X
\in \f h$), which allows us to identify $\f h$ with $\fhd$. Under
this identification, we have $\iota(H_\al) = \frac{2 \al}{(\al,
\al)} =: \al^{\vee}$, where $H_\al \in \fh$ is the simple coroot corresponding to $\al \in
\Pi$.

For each $\al \in \Phi$, set $s_\al(\xi) = \xi - (\xi, \al^\vee)
\al$ (for $\xi \in \fhd$) to be the reflection through the
hyperplane defined by $\al^\vee$.  For $\al_i \in \Pi$, let $s_i :=
s_{\al_i}$, be the simple reflection defined by $\al_i$.
Define $\Pi_\fk$ to be the set of simple roots in $\Phi^+_c$ and let $W_\fk$ denote
the Weyl group generated the reflections defined by $\Pi_\fk$.
Let $W_c = \langle s_\al | \al\in \Pi_c \rangle$ be the parabolic
subgroup of $W_\fg$ defined by the compact simple $\fg$-roots.  Note that
$W_c \subseteq W_\fk$, but we do not have equality in general.
For $w \in W_\fk$, set $\ell(w) := |w(\Phi^+_c) \cap \Phi^-_c|$.  Note that there is also a
length function on $W_\fg$ (denoted $\ell_\fg$) but $\ell$ refers to $W_\fk$.

A weight $\xi \in \fhd$ is said to be $\fk$-dominant (resp.
$\fg$-dominant) if $(\xi, \al) \geq 0$ for all $\al \in \Pi_\fk$ (resp
$\al \in \Pi$).  A weight $\xi\in \fhd$ is $\fg$-regular (resp. $\fk$-regular)
if $(\xi,\al)\neq 0$ for all $\al \in \Phi$ (resp. $\al \in \Phi_c$).
The integral weight lattice for $\fg$ is denoted by the set \\
$P(\fg) = \{\xi \in \fhd| (\xi, \al^{\vee}) \in\bbZ \mbox{ for all }
\al \in \Pi_\fg \}.$ Similarly we let $P(\fk)$ denote the abelian
group of integral weights for $\fk$ corresponding to $\Pi_\fk$. Let
the set of $\fk$- and $\fg$-dominant integral weights be denoted by
$P_+(\fk)$ and $P_+(\fg)$ respectively.  To each element $\delta \in
P_+(\fk)$ (resp. $P_+(\fg)$), let $L_\fk(\de)$ (resp. $L_\fg(\de)$)
denote the finite dimensional representation of $\fk$ (resp. $\fg$)
with highest weight $\de$.

Next, let $Q : P(\fk) \rightarrow \bbN$ denote the $\Phi^+_{nc}$-partition
function.  That is, if $\xi \in P(\fk)$ then $Q(\xi)$ is the number
of ways of writing $\xi$ as a sum of noncompact positive roots.  Put other way:
there exists an algebraic torus, $T$, such that to each $\mu \in P(\fk)$
there corresponds a linear character of $T$, denoted $e^\mu$, with
differential $\mu$.  Thus, $Q$ defines the coefficients of the product:
\[
    \sum_{\xi \in \fhd} Q(\xi) e^\xi = \prod_{\ga \in \Phi^+_{nc}}
    (1-e^\xi)^{-1}.
\]

Finally, we define the Blattner formula.  For $\de, \mu \in P(\fk)$,
\[
B(\de,\mu) := \sum_{w \in W_\fk} (-1)^{\ell(w)} Q(
w(\de+\rho)-\rho-\mu).
\]

It is convenient to introduce the notation $w.\xi =
w(\xi+\rho)-\rho$ for $w \in W_\fk$ and $\xi \in \fhd$.  It is easy
to see that $B(v.\de,\mu) = (-1)^{\ell(v)} B(\de,\mu)$.  Since for
all $\de \in P(\fk)$ there exists $v \in W_\fk$ such that $v.\de \in
P_+(\fk)$, we will assume that $\de \in P_+(\fk)$.

Historically, Blattner's formula arises out of the study of the
discrete series and its generalizations (see \cites{En, EW, EV, HC, HS, GrWa, Sc, Vo, Wa1, Wa2}).
\begin{thm}[See \cite{HS}]\label{thm_HS}
Assume $\la = \la(\mu) := \mu - \rho_{nc} + \rho_c$ is
$\fg$-dominant and $\fg$-regular. Then, $B(\de,\mu)$ is the
multiplicity of the finite dimensional $\fk$-representation,
$L_\fk(\de)$, in the discrete series representation of $G$ with
Harish-Chandra parameter $\la$.
\end{thm}

In this paper, we do not impose the $\fg$-dominant regular condition
on $\la(\mu)$. This is natural from the point of view of representation
theory as it is related to the coherent continuation of the discrete
series (see \cite{Sc}, \cite{Vo} and \cite{Zu}).

\bigskip
From our point of view, the goal is to understand the Blattner formula in as combinatorial fashion as possible.
Thus it is convenient to introduce the following generating function:

\newpage
\begin{dfn*}  For $\de \in P_+(\fk)$ we define the formal series:
\[ \bl(\delta) := \sum_{\mu \in \fhd} B(\delta,\mu) e^\mu. \]
\end{dfn*}

\noindent The main result of this paper is Proposition \ref{prop_main} of Section \ref{sec_main}, which states:
For $\de \in P_+(\fk)$,
\[  \bl(\delta) = \ch L_\fk(\delta) \, \frac{\prod_{\ga \in \Phi^+_{c}} (1-e^{-\ga})}{\prod_{\ga \in \Phi^+_{nc}}
(1-e^{-\ga})}, \]
where $\ch L_\fk(\de)$ denotes the character of $L_\fk(\de)$.
\bigskip

Of particular interest are the cases where $\Pi_c \neq \emptyset$,
which we address in Section \ref{sec_sym}.  From the point of view
of representation theory these include, for example, the holomorphic
and Borel-de Siebenthal discrete series (see \cite{HS}). More
recently, the latter has been addressed in \cite{GrWa}.

The Blattner formula for the case of $\Pi_c = \emptyset$ is often
particularly difficult to compute explicitly when compared to, say,
the cases corresponding to holomorphic discrete series. The $\Pi_c =
\emptyset$ case corresponds to the \emph{generic} discrete series of
the corresponding real semisimple Lie group.  In Section
\ref{sec_G2} we explore this situation in some detail for the Lie
algebra $\rm G_2$.

\bigskip

Finally, in light of Theorem \ref{thm_HS} one may observe that if
$\de \in P_+(\fk)$ and $\la(\mu)$ is $\fg$-dominant regular then
$B(\de,\mu) \geq 0$.  Our goal is to investigate the positivity
of Blattner's formula using combinatorial methods.  Of particular
interest is the positivity when we relax the $\fg$-dominance
condition on $\la(\mu)$.  Some results in this direction are
suggested by the recent work of Penkov and the second author
(see \cite{PZ}).

In Section \ref{sec_sym} we prove the existence of a skew symmetry
of Blattner's formula that exists whenever $\Pi_c \neq \emptyset$.
Thus, the condition that $\Pi_c = \emptyset$ is necessary for
$B(\de,\mu) \geq 0$ for all $\de \in P_+(\fk)$ and $\mu \in P(\fk)$.
In the situation where $\Pi_c =\emptyset$ we introduce the
following:

\begin{dfn*}  We say that a semisimple Lie algebra is $\bl$-positive
if the Blattner formula corresponding to the $\bbZ_2$-gradation with
$\Pi_c = \emptyset$ has the property that:
\[ B(\de,\mu)\geq 0 \hspace{1cm} \mbox{for all $\de \in P_+(\fk)$ and $\mu \in P(\fk)$.} \]
\end{dfn*}

The terminology stems from the fact that a simple Lie algebra is
$\bl$-positive if and only if the coefficients of $\bl(\de)$ are
non-negative for all $\de \in P_+(\fk)$.  Since the character of
$L_\fk(\de)$ can be written as a non-negative integer combination of
characters of $T$, we have that $\bl(\de)$ has non-negative integer
coefficients if $\bl(0)$ does.  Thus the question of
$\bl$-positivity reduces to the case of $\de = 0$.  In Section \ref{sec_positive},
it is shown that the only $\bl$-positive
simple Lie algebras are of type $\rm A_1$, $\rm A_2$, $\rm A_3$, $\rm B_2$,
$\rm C_3$, $\rm D_4$ and $\rm G_2$.  We prove this result by
examining the coefficients of $\bl(0)$.

\newpage
\section{Proof of the Main Result}\label{sec_main}

\begin{prop}\label{prop_main}  For $\de \in P_+(\fk)$,
\[  \bl(\delta) = \ch \; L_\fk(\delta) \, \frac{\prod_{\ga \in \Phi^+_{c}} (1-e^{-\ga})}{\prod_{\ga \in \Phi^+_{nc}}
(1-e^{-\ga})}, \]
where $\ch \; L_\fk(\de)$ denotes the character of $L_\fk(\de)$.
\end{prop}
\begin{proof}
From the definition of Blattner's formula we have:
\[
\bl(\delta) = \sum_{\mu \in \fhd} \sum_{w \in W_\fk} (-1)^{\ell(w)} Q(
w(\delta+\rho)-\rho-\mu) e^\mu.
\]
First we make the substitution, $\mu = w(\delta+\rho)-\rho-\xi$, and reorganize the sum:
\begin{eqnarray*}
\bl(\delta)
& = & \sum_{w \in W_\fk} (-1)^{\ell(w)} \sum_{\mu \in \fhd} Q(\xi) e^{w(\delta+\rho)-\rho-\xi} \\
& = & \sum_{w \in W_\fk} (-1)^{\ell(w)} e^{w(\delta+\rho)-\rho} \sum_{\mu \in \fhd} Q(\xi) e^{-\xi} \\
& = & \frac{\sum_{w \in W_\fk} (-1)^{\ell(w)} e^{w(\delta+\rho)-\rho}}{ \prod_{\ga \in \Phi^+_{nc}} (1-e^{-\ga})} \\
& = & \frac{\sum_{w \in W_\fk} (-1)^{\ell(w)} e^{w(\delta+\rho)-\rho}} {\prod_{\ga \in \Phi^+_{c}} (1-e^{-\ga})} \frac{\prod_{\ga \in \Phi^+_{c}} (1-e^{-\ga})}{ \prod_{\ga \in \Phi^+_{nc}} (1-e^{-\ga})}.
\end{eqnarray*}

As is well known, the character may be expressed using Weyl's
formula (see \cite{GW, Kn}) as in the following:
\[
    \ch \; L_\fk(\de) = \frac{\sum_{w \in W_\fk} (-1)^{\ell(w)} e^{w(\de + \rho)-\rho}}{\prod_{\ga \in \Phi^+_c} (1 - e^{-\ga})}.
\]
\end{proof}

\bigskip

The result allows us to compute the Blattner formula as follows:
\[ \sum_{\mu \in \fhd} B(\de,\mu) e^\mu = \ch L_\fk(\de) \, \sum_{\nu \in \fhd} B(0,\nu)e^\nu. \]
Let $\ch L_\fk(\de) = \sum_{\ga \in \fhd} m_\ga e^\ga$ and we obtain:
\[ \sum_{\mu \in \fhd} B(\de,\mu) e^\mu = \sum_{\ga, \nu \in \fhd} m_\ga B(0,\nu) e^{\ga+\nu}. \]
Thus for $\de \in P_+(\fk)$ and $\mu \in P(\fk)$ we have
\[ B(\de,\mu)  = \sum_{\ga \in P(\fk)} m_\ga B(0,\mu-\ga). \]
Note that the numbers $m_\ga$ are weight multiplicities for the representation $L_\fk(\delta)$.

\section{(Skew-)Symmetries of Blattner's Formula}\label{sec_sym}

The main result of this section is:
\begin{prop}\label{thm_sym} For $v \in W_c$:
\[
    \begin{array}{ll}
        B(\delta, \mu) = B(\delta,v. \mu) & \text{if $\ell(v)$ is even}, \\
        B(\delta, \mu) =-B(\delta,v. \mu) & \text{if $\ell(v)$ is odd}.
    \end{array}
\]
\end{prop}

Although this is well known to experts, we include our proof as it requires very little technical machinery.

\begin{dfn*}  For $w \in W_\fg$ and $\xi \in \fhd$ let
$Q_w(\xi) := Q(w^{-1}\xi)$.
\end{dfn*}

\begin{lemma}\label{lemma_Q} If $w(\Phi^+_{nc}) = \Phi^+_{nc}$ then $Q_w =
Q$.
\end{lemma}
\begin{proof}  It is enough to show
\[ \sum_{\xi\in \fhd} Q_w(\xi) e^\xi = \sum_{\xi\in \fhd} Q(\xi) e^\xi, \]
which follows from the following calculation:
\[\begin{array}{lllll}
    \sum_{\xi} Q_w(\xi) e^\xi &=& \sum_\xi Q(\xi) e^{w(\xi)} &=& w\left( \sum_{\xi} Q(\xi) e^\xi \right) \\
                              &=& w\left( \prod_{\al \in \Phi^+_{nc}} \left(1-e^\al\right)^{-1} \right)
                              &=& \prod_{\al \in \Phi^+_{nc}} \left(1-e^{w(\al)}\right)^{-1} \\
                              &=& \prod_{w^{-1}(\al) \in \Phi^+_{nc}} \left(1-e^\al\right)^{-1}
                              &=& \prod_{\al \in w(\Phi^+_{nc})} \left(1-e^\al \right)^{-1}.
\end{array}\]
\end{proof}

\begin{lemma}\label{lemma_roots} For all $w \in W_c$, $w(\Phi^+_{nc}) =
\Phi^+_{nc}$.
\end{lemma}
\begin{proof}
Note that $W_c$ is the Weyl group of a reductive Levi factor
$\mathfrak l$ of a parabolic subalgebra $\mathfrak q \subseteq
\mathfrak g$. We have a generalized triangular decomposition
$\mathfrak g = \mathfrak u^- \oplus \mathfrak l \oplus \mathfrak
u^+$ with $\mathfrak q=\mathfrak l \oplus \mathfrak u^+$. The
noncompact root spaces contained in $\mathfrak u^+$ are positive.
Furthermore, all noncompact positive roots spaces are contained in
$\mathfrak u^+$ because $\mathfrak l \subseteq \fk$. The Lie algebra
$u^+$ is an $\mathfrak l$-module and therefore the weights are
preserved by $W_c$.  It is this fact that implies that $W_c$ takes
positive noncompact roots to positive roots.  We now need to show that
$W_c$ takes noncompact roots to noncompact roots.

For roots $\beta \in \Phi_{nc}$ and $\al \in \Phi_c$, we have the
formula
\[
    s_\al(\beta) = \beta - \frac{2(\al,\beta)}{(\al,\al)} \al
\] with $\frac{2(\al,\beta)}{(\al,\al)}$ an integer.  Thus,
the reflection of a noncompact root across a hyperplane defined by
a compact root is noncompact.  The reflections generate $W_c$.
\end{proof}

We now prove the main proposition of this section.

\begin{proof}[Proof of Proposition \ref{thm_sym}]
Let $v \in W_c$:
\begin{eqnarray*}
B(\delta, \mu) &=& \sum_{w\in W_{\fk}} (-1)^{\ell(w)} Q( w.\delta - \mu) \\
&=& \sum_{ w \in v(W_{\fk})} (-1)^{\ell(v^{-1}w)} Q( (v^{-1}w).
\delta - \mu).
\end{eqnarray*}
By definition, $\ell(v) = \ell(v^{-1})$.  Combining this with the
definition of the ``dot" action we obtain:

\begin{eqnarray*}
B(\delta,\mu) &=& (-1)^{\ell(v)} \sum_{ w \in W_{\fk}}
(-1)^{\ell(w)} Q( v^{-1}w(\delta +
\rho) - \rho - \mu) \\
&=& (-1)^{\ell(v)} \sum_{ w \in W_{\fk}} (-1)^{\ell(w)} Q\left(
v\left(v^{-1}w(\delta + \rho) - \rho - \mu \right) \right).
\end{eqnarray*}
Next, we use the fact that $Q = Q_{v^{-1}}$ using Lemmas \ref{lemma_Q} and \ref{lemma_roots}.
The rest is a calculation.
\begin{eqnarray*}
B(\delta,\mu) &=& (-1)^{\ell(v)} \sum_{ w \in W_{\fk}}
(-1)^{\ell(w)} Q\left(w(\delta
+ \rho) - v(\mu + \rho) \right) \\
&=& (-1)^{\ell(v)} \sum_{ w \in W_{\fk}} (-1)^{\ell(w)}
Q\left(w(\delta + \rho) - \rho - v(\mu
+ \rho) +\rho  \right) \\
&=& (-1)^{\ell(v)} \sum_{ w \in W_{\fk}} (-1)^{\ell(w)}
Q\left(w(\delta + \rho) - \rho - \left(v(\mu + \rho) -\rho\right)
\right)
\\ &=& (-1)^{\ell(v)} \sum_{ w \in W_{\fk}} (-1)^{\ell(w)}
Q\left(w . \delta - v . \mu \right) = (-1)^{\ell(v)} B(\delta, v.
\mu).
\end{eqnarray*}
\end{proof}

\section{The case of $\rm G_2$}\label{sec_G2}

\bigskip
The following is a complete calculation of $\bl(0)$ for the case of
the Lie algebra $\fg := \rm G_2$ when $\Pi_c = \emptyset$. Let $\al$
and $\be$ be a choice of noncompact simple roots for $\rm G_2$ with
$\al$ long and $\be$ short. The compact positive roots are $\Phi_c^+
= \{\al+\be, \al+3\be \}$, while the noncompact positive roots are
$\Phi_{nc}^+ = \{\al,\be,\al+2\be,2\al+3\be\}$. Denote the
$\Phi_{nc}^+$-partition function by $Q : \fhd \rightarrow \bbZ$. We
have:
\[
   \q := \sum_{\xi \in \fhd} Q(\xi) e^{-\xi} = \prod_{\ga \in \Phi_{nc}^+}
    (1-e^{-\ga})^{-1}.
\]

Let $x = e^{-\al}$ and $y = e^{-\be}$.  Thus,
\[
\q = \frac{1}{(1-x)(1-y)(1-x y^2)(1-x^2 y^3)}.
\]
Let the simple reflection corresponding to $\al$ (resp. $\be$) be
$s_1$ (resp. $s_2$).  We have four terms in Blattner's formula for
$\delta = 0$:
\[
B(0,\mu) = Q(-\mu) - Q(s_1\rho-\rho-\mu) - Q(s_2\rho-\rho-\mu) +
Q(s_1 s_2\rho - \rho - \mu),
\]
and our goal will be to close the sum $\bl(0) := \sum_{\mu \in \fhd}
B(0,\mu) e^{\mu}$.  We will do this by multiplying by $e^{\mu}$
and summing over $\mu$ for each of the four terms. Observe that:
\[
\sum_{\xi \in \fhd} Q(\xi) e^{-\xi} = \sum_{\xi \in \fhd} Q(-\xi)
e^{\xi}.
\]
Thus, $\q = \sum_{\xi \in \fhd} Q(-\xi) e^{\xi}$. Next consider the
sum:
\[
    T_1 := \sum_{\mu \in \fhd} Q(s_1\rho-\rho-\xi) e^{\xi}.
\]  We make the substitution $-\mu = s_1\rho - \rho - \xi$ so that the
above sum becomes:
\[
    \sum_{\mu \in \fhd} Q(-\mu)e^{s_1\rho-\rho+\mu} =
    e^{s_1 \rho - \rho} \sum_{\mu \in \fhd} Q(-\mu) e^{\mu}.
\]
Thus the above sum is equal to $e^{s_1 \rho - \rho} \q$, which we
denote by ${\bf T_1}$.  Similarly, we set ${\bf T_3} := e^{s_2 \rho
- \rho} \q$ and ${\bf T_3} := e^{s_1 s_2 \rho - \rho} \q$, and
${\bf T_0} := \q$. Thus,
$
\bl(0) = {\bf T_0} - {\bf T_1} - {\bf T_2} + {\bf T_3}.
$

Now we write the above in terms of $x$ and $y$.  Note that $\rho =
\al + 2 \be$, $s_1 \rho = \be$, $s_2 \rho = -\be$, and from these we
can easily see:
\begin{eqnarray*}
    e^{s_1     \rho - \rho} & = x y,    \\
    e^{s_2     \rho - \rho} & = x y^3,  \\
    e^{s_1 s_2 \rho - \rho} & = x^2 y^4 .
\end{eqnarray*}

Putting everything together we see $\bl(0) = (1- x y)(1- x y^3) \q$. Or
equivalently, if $\mu = -k \al - \ell \be$  then the value of
$B(0,\mu)$ is the coefficient of $x^k y^\ell$ in:
\begin{eqnarray*}
\bl(0) &=& \frac{(1- x y)(1- x y^3)}{(1-x)(1-y)(1-x y^2)(1-x^2 y^3)} \\
       &=&  \frac{1}{(1-x)(1-x y^2)} + \frac{y}{(1-y)(1-x^2 y^3)}.
\end{eqnarray*}
Note that in the latter expression, it is clear that the
coefficients are the series expansion are positive. The positivity
of the coefficients of $\bl(\delta)$ follows from the positivity of
the coefficients of $\bl(0)$.  The question of positivity for a
general semisimple Lie algebra will be addressed in Section
\ref{sec_positive}.

It is important to note that as we change $\Pi_c$ the $\bl(0)$ changes as well.  For
example, when $\Pi_c=\{ \be \}$,
\begin{eqnarray*}
\bl(0) &=& \frac{(1-x^2 y^3)(1-y)}{(1-x)(1-xy)(1-x y^2)(1-xy^3)} \\
&=&
\frac{1}{(1-x)(1-xy^2)}+\frac{x}{(1-x^2)(1-xy)}-\frac{x+y+xy^2}{(1-x^2)(1-xy^3)}.
\end{eqnarray*}
which we can easily see does \emph{not} have non-negative integer coefficients.
We leave it as an exercise to the reader to compute $\bl(0)$ when $\Pi_c = \{ \al \}$.

The generating function for other $\delta$ involve multiplying the
above product by a polynomial in $x^{\pm \frac{1}{2}}, y^{\pm \frac{1}{2}}$ that represents
the character of the corresponding irreducible finite dimensional representation of $\fk = \fso_4$.

\bigskip

One can plot the coefficients of our formal series, as we do next.
In all pictures we have labeled the scale on the axes and normalized
the short root ($\be$) to have length 1 and be positioned at 3 o'clock.
The long root $\al$ is at 10 o'clock.

The first two images (Figures 1 and 2)
are for the case with $\Pi_c = \emptyset$ (generic) and correspond to
$\de = 0$ and  $\de$ the highest weight of the standard $4$-dimensional
representation of $\fso_4$.  We also display the same two $\de$'s in the case
when $\Pi_c = \{ \be \}$ (Borel-de Siebenthal) in Figures 3 and 4.
These latter two figures clearly display a skew-symmetry
addressed Section \ref{sec_sym}.

\newpage

\begin{center}
\includegraphics[width=.6\textwidth]{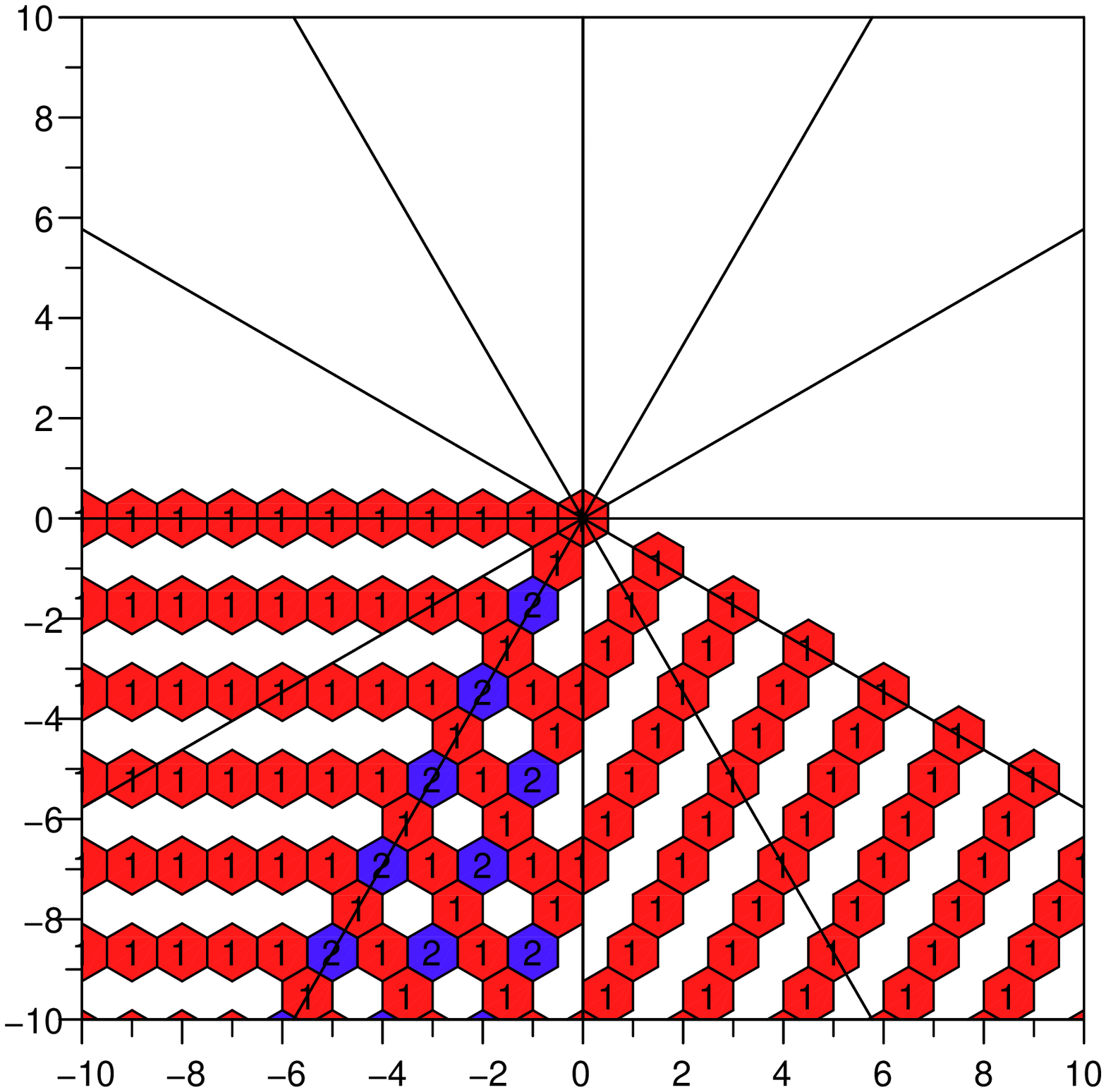}
\center{Figure 1: Coefficients of $\bl(0)$ when $\Pi_c =\emptyset$.}

\vspace{2cm}

\includegraphics[width=.6\textwidth]{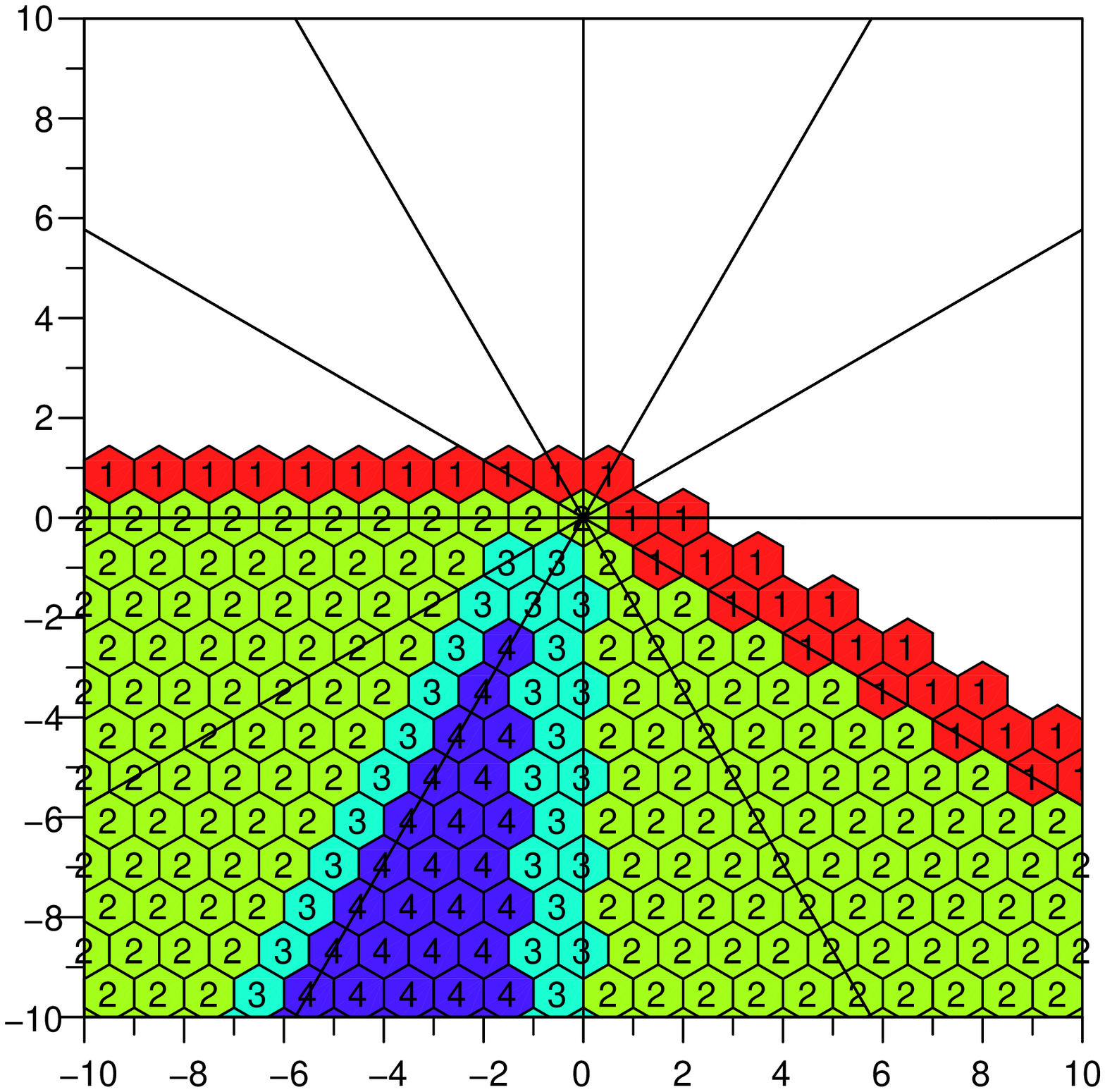}
\center{Figure 2: $\bl(\de)$ for $L_{\fso_4}(\de)$ the standard rep. ($\Pi_c = \emptyset$)}
\end{center}

\newpage

\begin{center}
\includegraphics[width=.6\textwidth]{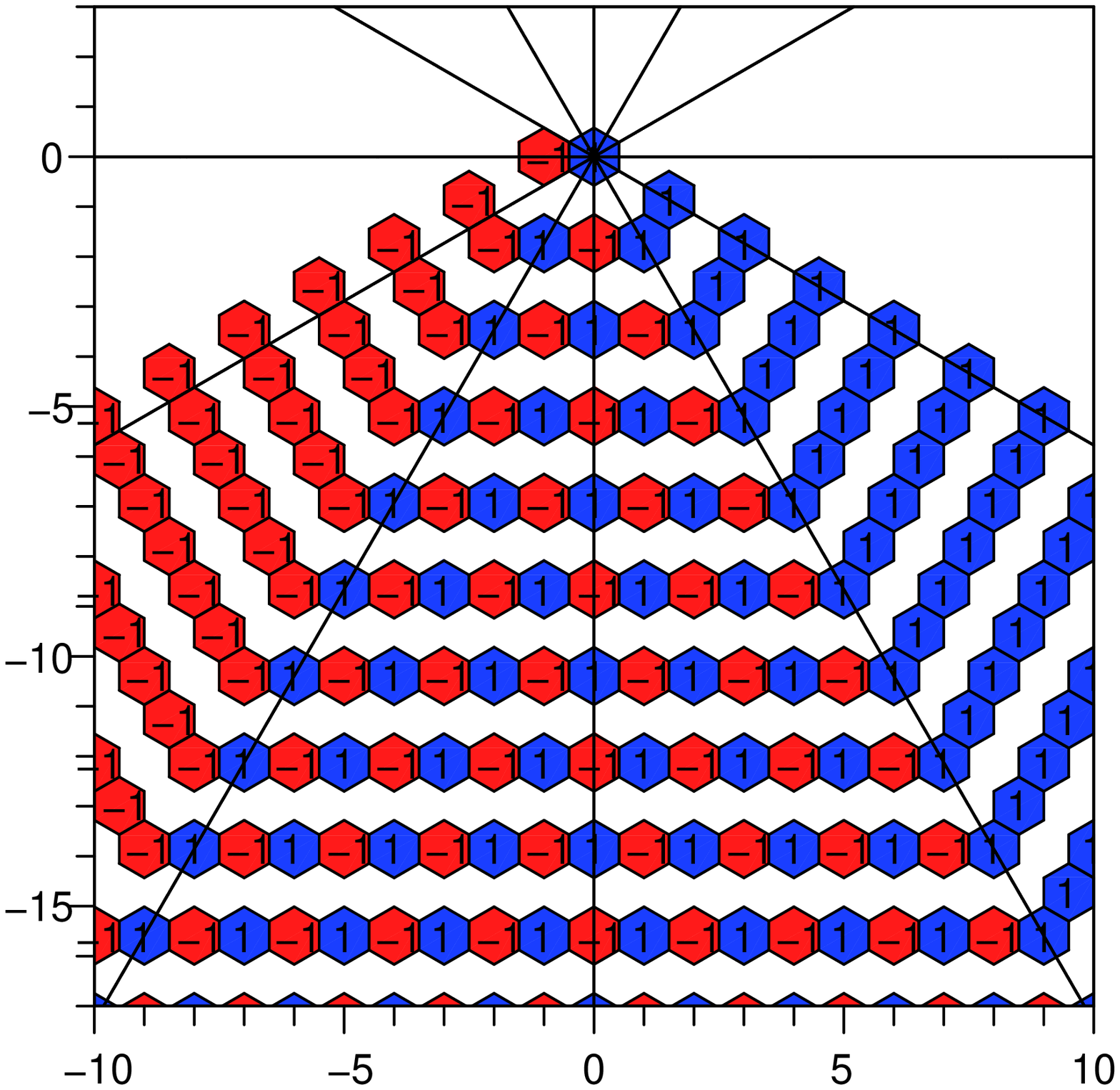}
\center{Figure 3: Coefficients of $\bl(0)$ when $\Pi_c = \{ \be\}$.}

\vspace{2cm}

\includegraphics[width=.6\textwidth]{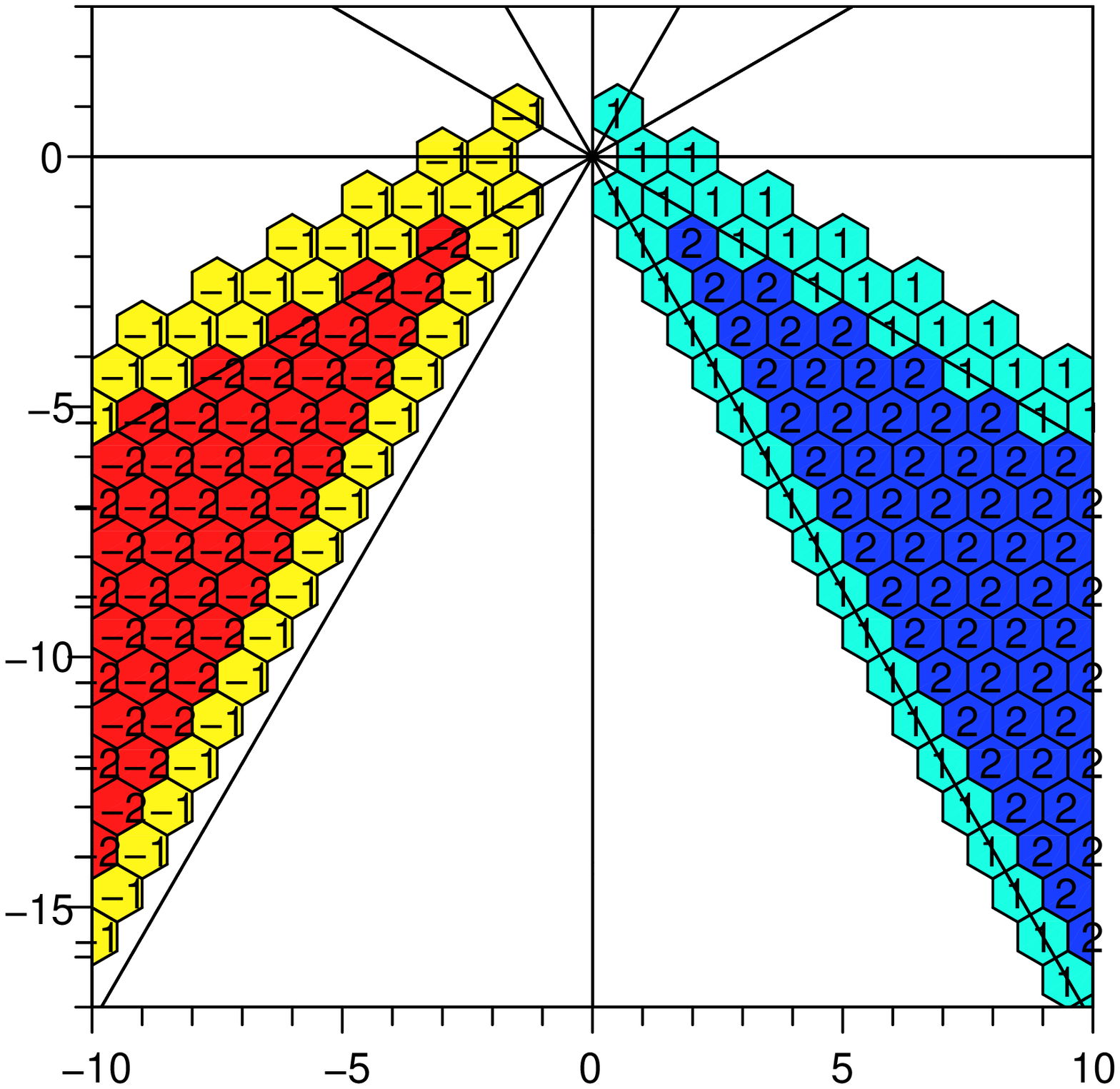}
\center{Figure 4: $\bl(\de)$ for $L_{\fso_4}(\de)$ the standard rep. ($\Pi_c = \{ \be \}$)}
\end{center}

\newpage

\section{$\bl$-positive simple Lie algebras}\label{sec_positive}

We have the following:
\begin{prop}\label{prop_seven}  The only $\bl$-positive simple Lie algebras are
of type $\rm A_1$, $\rm A_2$, $\rm A_3$, $\rm B_2$, $\rm C_3$, $\rm D_4$ and
$\rm G_2$.
\end{prop}

To prove the proposition we use the following lemma to exclude the
other cases.
\begin{lemma}\label{lem_levi}  Let $\fl$ be the semisimple Levi factor of a parabolic
subalgebra of $\fg$.  If $\fg$ is $\bl$-positive then $\fl$ is also
$\bl$-positive.
\end{lemma}
\begin{proof}
It is a consequence of the main proposition that for any rank $r$
simple Lie algebra $\fg$ we may express $\bl(0)$ as a rational
function in $x_1,\cdots,x_r$ where $x_i = e^{-\al_i}$ (recall that
$\al_i$ is the $i^{\rm{th}}$-simple root for $\fg$). For any $j$, if
we set $x_j = 0$ in $\bl(0)$ then the resulting expression is
$\bl(0)$ for a semisimple subalgebra of $\fg$.  It is not hard to
see that this subalgebra is the Levi factor of the maximal parabolic
of $\fg$ corresponding to $\al_j$.  More generally, let $S \subseteq
\Pi$. If we set $x_i = 0$ in $\bl(0)$ for all $\al_i \in S$ then the
resulting expression is $\bl(0)$ for the Levi subalgebra, $\fl_S$,
of the corresponding parabolic.  Note that the terms in the series
expansion for $\bl(0)$ for $\fl_S$ are also terms in the series
expansion of $\bl(0)$ for $\fg$.  Thus a negative coefficient in the
former implies a negative coefficient in the latter.
\end{proof}

The proof of Proposition \ref{prop_seven} involves a cases-by-case
analysis using the main proposition expressing $\bl(0)$ as a
product, Lemma \ref{lem_levi} and the computer algebra package,
MAPLE.    Recall that the $\bl$-positivity of $\rm G_2$ was proved
in Section \ref{sec_G2}.

\subsection{Type A}

For $\fsl_2$, we have: $ \bl(0) = \frac {1}{1 - {x_{1}}}. $ For
$\fsl_3$:
\[
\bl(0)=\frac {1 - {x_{1}}\,{x_{2}}}{(1 - {x_{1}})\,(1 - {x_{2}})} =
1 + \frac {{x_{1}}}{1 - {x_{1}}} + \frac {{x_{2}}}{1 - {x_{2}}}.
\]
And for $\fsl_4$:
\[
\bl(0) = \frac {(1 - {x_{1}}\,{x_{2}})\,(1
 - {x_{2}}\,{x_{3}})}{(1 - {x_{1}})\,(1 - {x_{2}})\,(1 - {x_{3}})
\,(1 - {x_{1}}\,{x_{2}}\,{x_{3}})} = \frac {1}{(1 - {x_{3}})\,(1 -
{x _{1}})}  + \frac {{x_{2}}}{(1 - {x_{2}})\,(1 - {x
_{1}}\,{x_{2}}\,{x_{3}})}.
\]
The partial fraction on the right of each of these examples
establishes that the coefficients are indeed positive.  Next we
consider $\fsl_5$ where we have a negative result. First we note:
\[
\bl(0) = \frac {(1 - {x_{1}}\,{x_{2}})\,(1
 - {x_{2}}\,{x_{3}})\,(1 - {x_{3}}\,{x_{4}})\,(1 - {x_{1}}\,{x_{2
}}\,{x_{3}}\,{x_{4}})}{(1 - {x_{1}})\,(1 - {x_{2}})\,(1 - {x_{3}}
)\,(1 - {x_{4}} )\,(1 - {x_{1}}\,{x_{2}}\,{x_{3}})\,(1 -
{x_{2}}\,{x_{3}}\,{x_{4} })},
\]
for $\fsl_5$.  We expand in a formal power series and observe that
the coefficient of $x_1 x_2^2 x_3^2 x_4$ is $-1$.  This means that
$\bl$-positivity fails for this example.  We then also see failure
of $\bl$-positivity for any simple Lie algebra which has $\rm A_4$
as a Levi factor of a parabolic. Thus, we exclude all higher rank
type A examples as well as $\rm B_n$ ($n\geq5$), $\rm C_n$ ($n \geq
5$), $\rm D_n$ ($n \geq 5$), $\rm E_6$, $\rm E_7$ and $\rm E_8$.

\subsection{Type B}
We only need to examine $\fso_5$ and $\fso_7$. For $\fso_5$ we have:
\[
{\displaystyle \bl(0) = \frac {1 - {x_{1}}\,{x_{2}}}{(1 - {
x_{1}})\,(1 - {x_{2}})\,(1 - {x_{1}}\,{x_{2}}^{2})}} =
{\displaystyle \frac {1}{( 1 - {x_{2}}^{2})\,(1
 - {x_{1}})}}  + {\displaystyle \frac {{x_{2}}}{(1 - {x_{2}}^{2}
)\,(1 - {x_{1}}\,{x_{2}}^{2})}}.
\]
Thus the coefficients of the series expansion are non-negative.

Now consider $\fso_7$:
\[
\bl(0) = \frac {(1 - {x_{1}}\,{x_{2}})\,(1 - {x_{2}}\,{x_{3 }})\,(1
- {x_{1}}\,{x_{2}}\,{x_{3}}^{2})}{(1 - {x_{1}})\,(1 - {x _{2}})\,(1
- {x_{3}})\,(1 - {x_{1}}\,{x_{2}}\,{x_{3}})\,(1 - {x_{
2}}\,{x_{3}}^{2})\,(1 - {x_{1}}\,{x_{2}}^{2}\,{x_{3}}^{2})}.
\]
Upon expansion we see that the coefficient of $x_1^2 x_2^3 x_3^3$ is
$-1$.  Thus, we may exclude this and higher rank type $\rm B$
examples as they have $B_3$ as the Levi factor of a parabolic.  In
particular, we may exclude $B_4$, which is the only type $\rm B$
example that has not been excluded yet.  We also may exclude $\rm
F_4$ for the same reason.

\subsection{Type C}  We must examine $\bl(0)$ for $\fsp_6$ and
$\fsp_8$ as these are the only examples not yet addressed.  For
$\fsp_{6}$ (ie: $C_3$), we have:
\[
\bl(0)= \frac {(1 - {x_{1}}\,{x_{2}})\,(1 - {x_{2}}\,{x_{3 }})\,(1 -
{x_{1}}\,{x_{2}}^{2}\,{x_{3}})}{(1 - {x_{1}})\,(1 - {x _{2}})\,(1 -
{x_{3}})\,(1 - {x_{1}}\,{x_{2}}\,{x_{3}})\,(1 - {x_{
2}}^{2}\,{x_{3}})\,(1 - {x_{1}}^{2}\,{x_{2}}^{2}\,{x_{3}})}.
\]
The coefficients are positive as the above expression is equal to:
\begin{eqnarray*}
&\frac {1}{(1 - {x_{1}})\,(1 - {x_{2}}^{2}\,{x_{3}})} +
 \frac {{x_{1}}\,{x_{3}}}{(1 - {x_{1}})\,(1 - {x_{3
}})\,(1 - {x_{1}}^{2}\,{x_{2}}^{2}\,{x_{3}})} + \frac
{{x_{1}}^{2}\,{x_{2}} ^{3}\,{x_{3}}^{2}}{(1 -
{x_{2}}^{2}\,{x_{3}})\,(1 - {x_{1}}^{2}\,{x_{2}}^{2}\,{x_{3}})\,(1 -
{x_{1}}\,{x_{ 2}}\,{x_{3}})}+\\ & \frac {{x_{3}}}{(1 -
{x_{2}}^{2}\,{x_{3}})\, (1 - {x_{1}}^{2}\,{x_{2}}^{2}\,{x_{3}})\,(1
- {x_{3}})} + \frac {{x_{2}}}{(1 - {x_{2}}^{2}\,{x_{3}})\,(1 -
{x_{1}}^{2}\,{x_{2}}^{2}\,{x_{3}})\,(1 - {x_{2}})}.
\end{eqnarray*}

However, we do not have $\bl$-positivity for $\fsp_8$, as we see
$-1$ as the coefficient of $x_1 x_2^3 x_3^3 x_4^2$ in $\bl(0) =
\frac{NUM}{DEN}$ where:
\begin{eqnarray*}
& NUM  := & (1 - {x_{1}}\,{x_{2}})\,(1 - {x_{2}}\,{x_{3}})\,(1 -
{x_{3}}\,{x _{4}})\, (1 - {x_{1}}\,{x_{2}}\,{x_{3}}\,{x_{4}})  \\
& & (1 - {x_{2}}\,{ x_{3}}^{2}\,{x_{4}})\,(1 -
{x_{1}}\,{x_{2}}^{2}\,{x_{3}}^{2}\,{x _{4}}), \; \mbox{and}
\end{eqnarray*}
\begin{eqnarray*}
& DEN  := & (1 - {x_{1}})\,(1 - {x_{2}})\,(1 - {x_{3}})\,(1 -
{x_{4}})\,(1 - {x_{1}}\,{x_{2}}\,{x_{3}})\,(1 -
{x_{2}}\,{x_{3}}\,{x_{4}})\,(1
 - {x_{3}}^{2}\,{x_{4}}) \\
& & (1 - {x_{1}}\,{x_{2}}\,{x_{3}}^{2}\,{x_{4}})\,(1 -
{x_{2}}^{2}\,{ x_{3}}^{2}\,{x_{4}})\,(1 -
{x_{1}}^{2}\,{x_{2}}^{2}\,{x_{3}}^{2} \,{x_{4}}).
\end{eqnarray*}

\subsection{Type D}  The only case left is $\fso_8$, where we have:
{\tiny
\begin{eqnarray*}
\bl(0)  &= & \frac {(1 - {x_{1}}\,{x_{2}})\,(1 - {x_{2}}\,{x_{3
}})\,(1 - {x_{2}}\,{x_{4}})\,(1 - {x_{1}}\,{x_{2}}\,{x_{3}}\,{x_{
4}})}{(1 - {x_{1}})\,(1 - {x_{2}})\,(1 - {x_{3}})\,(1 - {x_{4}})
\,(1 - {x_{1}}\,{x_{2}}\,{x_{3}})\,(1 - {x_{1}}\,{x_{2}}\,{x_{4}}
)\,(1 - {x_{2}}\,{x_{3}}\,{x_{4}})\,(1 - {x_{1}}\,{x_{2}}^{2}\,{x
_{3}}\,{x_{4}})} \\ & = & \frac {1}{(1 -
{x_{1}}\,{x_{2}}^{2}\,{x_{3}}\,{x_{ 4}})\,(1 - {x_{1}})\,(1 -
{x_{3}})\,(1 - {x_{4}})} + \frac {{x_{2}}}{(1 - {x_{1}}\,{x_{2}}\,{
x_{4}})\,(1 - {x_{2}})\,(1 - {x_{2}}\,{x_{3}}\,{x_{4}})\,(1 - {x
_{1}}\,{x_{2}}\,{x_{3}})}.
\end{eqnarray*}}
Thus the coefficients are non-negative integers.

\def\cprime{$'$} \def\cprime{$'$}


\begin{bibdiv}
\begin{biblist}


\bib{En}{article}{
   author={Enright, Thomas J.},
   title={On the fundamental series of a real semisimple Lie algebra: their
   irreducibility, resolutions and multiplicity formulae},
   journal={Ann. of Math. (2)},
   volume={110},
   date={1979},
   number={1},
   pages={1--82},
   issn={0003-486X},
   review={\MR{541329 (81a:17003)}},
}

\bib{EV}{article}{
   author={Enright, Thomas J.},
   author={Varadarajan, V. S.},
   title={On an infinitesimal characterization of the discrete series},
   journal={Ann. of Math. (2)},
   volume={102},
   date={1975},
   number={1},
   pages={1--15},
   issn={0003-486X},
   review={\MR{0476921 (57 \#16472)}},
}

\bib{EW}{article}{
   author={Enright, Thomas J.},
   author={Wallach, Nolan R.},
   title={The fundamental series of representations of a real semisimple Lie
   algebra},
   journal={Acta Math.},
   volume={140},
   date={1978},
   number={1-2},
   pages={1--32},
   issn={0001-5962},
   review={\MR{0476814 (57 \#16368)}},
}

\bib{GW}{book}{
    author={Goodman, R.},
    author={Wallach, N.R.},
     title={Representations and invariants of the classical groups},
 publisher={Cambridge University Press},
   address={Cambridge},
      date={1998},
      ISBN={0-521-58273-3},
    review={\MR{99b:20073}},
}

\bib{GrWa}{article}{
   author={Gross, B.},
   author={Wallach, N.},
   title={Restriction of small discrete series representations to symmetric
   subgroups},
   conference={
      title={The mathematical legacy of Harish-Chandra},
      address={Baltimore, MD},
      date={1998},
   },
   book={
      series={Proc. Sympos. Pure Math.},
      volume={68},
      publisher={Amer. Math. Soc.},
      place={Providence, RI},
   },
   date={2000},
   pages={255--272},
   review={\MR{1767899 (2001f:22042)}},
}


\bib{HC}{article}{
   author={Harish-Chandra},
   title={Discrete series for semisimple Lie groups. II. Explicit
   determination of the characters},
   journal={Acta Math.},
   volume={116},
   date={1966},
   pages={1--111},
   issn={0001-5962},
   review={\MR{0219666 (36 \#2745)}},
}

\bib{HS}{article}{
   author={Hecht, Henryk},
   author={Schmid, Wilfried},
   title={A proof of Blattner's conjecture},
   journal={Invent. Math.},
   volume={31},
   date={1975},
   number={2},
   pages={129--154},
   issn={0020-9910},
   review={\MR{0396855 (53 \#715)}},
}

\bib{Hu}{book}{
   author={Humphreys, James E.},
   title={Reflection groups and Coxeter groups},
   series={Cambridge Studies in Advanced Mathematics},
   volume={29},
   publisher={Cambridge University Press},
   place={Cambridge},
   date={1990},
   pages={xii+204},
   isbn={0-521-37510-X},
   review={\MR{1066460 (92h:20002)}},
}

\bib{Kn}{book}{
    author={Knapp, Anthony~W.},
     title={Lie groups beyond an introduction},
   edition={Second},
    series={Progress in Mathematics},
 publisher={Birkh\"auser Boston Inc.},
   address={Boston, MA},
      date={2002},
    volume={140},
      ISBN={0-8176-4259-5},
    review={\MR{2003c:22001}},
}

\bib{PZ}{article}{
   author={Penkov, Ivan},
   author={Zuckerman, Gregg},
   title={Generalized Harish-Chandra modules with generic minimal $\germ
   k$-type},
   journal={Asian J. Math.},
   volume={8},
   date={2004},
   number={4},
   pages={795--811},
   issn={1093-6106},
   review={\MR{2127949 (2005k:17007)}},
}

\bib{Sc}{article}{
   author={Schmid, Wilfried},
   title={$L\sp{2}$-cohomology and the discrete series},
   journal={Ann. of Math. (2)},
   volume={103},
   date={1976},
   number={2},
   pages={375--394},
   issn={0003-486X},
   review={\MR{0396856 (53 \#716)}},
}

\bib{Vo}{book}{
   author={Vogan, David A., Jr.},
   title={Representations of real reductive Lie groups},
   series={Progress in Mathematics},
   volume={15},
   publisher={Birkh\"auser Boston},
   place={Mass.},
   date={1981},
   pages={xvii+754},
   isbn={3-7643-3037-6},
   review={\MR{632407 (83c:22022)}},
}


\bib{Wa1}{book}{
   author={Wallach, Nolan R.},
   title={Real reductive groups. I},
   series={Pure and Applied Mathematics},
   volume={132},
   publisher={Academic Press Inc.},
   place={Boston, MA},
   date={1988},
   pages={xx+412},
   isbn={0-12-732960-9},
   review={\MR{929683 (89i:22029)}},
}

\bib{Wa2}{book}{
   author={Wallach, Nolan R.},
   title={Real reductive groups. II},
   series={Pure and Applied Mathematics},
   volume={132},
   publisher={Academic Press Inc.},
   place={Boston, MA},
   date={1992},
   pages={xiv+454},
   isbn={0-12-732961-7},
   review={\MR{1170566 (93m:22018)}},
}

\bib{Zu}{article}{
   author={Zuckerman, Gregg J.},
   title={Coherent translation of characters of semisimple Lie groups},
   conference={
      title={Proceedings of the International Congress of Mathematicians
      (Helsinki, 1978)},
   },
   book={
      publisher={Acad. Sci. Fennica},
      place={Helsinki},
   },
   date={1980},
   pages={721--724},
   review={\MR{562678 (81f:22025)}},
}

\end{biblist}
\end{bibdiv}
\end{document}